\documentclass[a4paper,12pt]{article}
\usepackage{amssymb} \usepackage{amsfonts}
\usepackage{amsmath,amsthm}
\usepackage{tikz,relsize}
\usepackage[cmtip,arrow]{xy}
\usepackage{pb-diagram,pb-xy}

\newcommand{\Z}{{\mathbb Z}}

\newcommand{\R}{{\mathbb R}}

\newtheorem{lm}{Lemma}[section]
\newtheorem{rem}{Remark}[section]
\newtheorem{theo}{Theorem}[section]

\theoremstyle{definition}
\newtheorem{defi}{Definition}[section]
\newtheorem{ex}{Example}[section]
\date{}

\begin{document}
\title{Spin-structures on real Bott manifolds}

\author{A. G\c{a}sior\footnote{Author is supported by the Polish National Science Center grant DEC-2013/09/B/ST1/04125}}


\maketitle

\hskip5mm
\section{Introduction}
Let $M^n$ be a flat manifold of dimension $n$, i.e. a compact connected Riemannian manifold without boundary
with zero sectional curvature. From the theorem of Bieberbach (\cite{Ch}, \cite{S3})
the fundamental group
$\pi_{1}(M^{n}) = \Gamma$ determines a short exact sequence:
\begin{equation}\label{ses}
0 \rightarrow \Z^{n} \rightarrow \Gamma \stackrel{p}\rightarrow
G \rightarrow 0,
\end{equation}
where
$\Z^{n}$ is a torsion free abelian group of rank $n$ and
$G$ is a finite group which
is isomorphic to the holonomy group of $M^{n}.$
The universal covering of $M^{n}$ is the Euclidean space $\R^{n}$
and hence $\Gamma$
is isomorphic to a discrete cocompact subgroup
of the isometry group $\operatorname{Isom}(\R^{n}) = \operatorname{O}(n)\times\R^n = E(n).$ In that case $p:\Gamma\to G$ is a projection on the first component of the semidirect product $O(n)\ltimes \mathbb R^n$ and $\pi_1(M_n)=\Gamma$ is a subfroup of $O(n)\ltimes \mathbb R^n$.
Conversely, given a short exact sequence of the form (\ref{ses}), it is known that
the group $\Gamma$ is (isomorphic to) the fundamental group of a flat manifold if and only if
$\Gamma$ is torsion free.
In this case $\Gamma$ is called a Bieberbach group.
We can define a holonomy representation $\phi:G\to \operatorname{GL}(n,\Z)$ by the formula:
\begin{equation}\label{holonomyrep}
\forall e\in\mathbb Z^n\;\forall g\in G,\phi(g)(e) = \tilde{g}e(\tilde{g})^{-1},
\end{equation}
where $p(\tilde{g})=g.$ In this article we shall consider Bieberbach groups of rank $n$ with holonomy group  $\Z_{2}^{k}$, $1\leq k\leq n-1$,
and  $\phi(\Z_{2}^{k})\subset D\subset \operatorname{GL}(n,\Z)$.
Here $D$ is the group of matrices with $\pm1$ on the diagonal.

Let
\begin{equation}\label{tower}
M_{n}\stackrel{\R P^1}\to M_{n-1}\stackrel{\R P^1}\to...\stackrel{\R P^1}\to M_{1}\stackrel{\R P^1}\to M_0 =
\{ \bullet\}
\end{equation}
be a sequence of real projective bundles such that $M_i\to M_{i-1}$, $i=1,2,\ldots,n$, is a projective
bundle of a Whitney sum of a real line bundle $L_{i-1}$ and the trivial line bundle over $M_{i-1}$.
The sequence (\ref{tower}) is called the real Bott tower and the top manifold $M_n$ is called the real
Bott manifold, \cite{CMO}.


Let $\gamma_i$ be the canonical line bundle over $M_i$ and we set $x_i = w_1(\gamma_i)$ ($w_1$ is the first Stiefel-Whitney class).
Since $H^1(M_{i-1},\Z_2)$ is additively generated by $x_1,x_2,..,x_{i-1}$ and $L_{i-1}$
is a line bundle over $M_{i-1},$ we can uniquely write
\begin{equation}\label{w1}
w_1(L_{i-1}) = \sum_{k=1}^{i-1} a_{ki}x_k
\end{equation}
where $a_{ki}\in \Z_2$ and $i = 2,3,...,n.$

From above we obtain the matrix $A = [a_{ki}]$ which is a $n\times n$ strictly upper triangular matrix
whose diagonal entries are $0$ and remaining entries are either
$0$ or $1.$
One can observe (see \cite{KM}) that the tower (\ref{tower}) is completly determined by the matrix $A$ and therefore we may denote the real Bott manifold $M_n$ by $M(A)$.
From \cite[Lemma 3.1]{KM} we can consider $M(A)$ as the orbit space  $M(A) = \R^n/\Gamma(A),$
where $\Gamma(A)\subset E(n)$ is generated by elements
\begin{equation}\label{gener}
s_{i} = \left(\left[
\begin{matrix}
1&0&0&.&.&...&0\\
0&1&0&.&.&...&0\\
.&.&.&.&.&...&\\
0&...&0&1&0&...&0\\
0&...&0&0&(-1)^{a_{i,i+1}}&...&0\\
.&.&.&.&.&...&\\
0&...&0&0&0&...&(-1)^{a_{i,n}}
\end{matrix}\right], \begin{pmatrix}
0\\
.\\
0\\
\frac{1}{2}\\
0\\
.\\
0\\
0
\end{pmatrix}\right)\in E(n),
\end{equation}
where $(-1)^{a_{i,i+1}}$ is in the $(i+1, i+1)$ position and $\frac{1}{2}$ is the $i-$th coordinate of the column, $i = 1,2,...,n-1.$
$s_{n} = \left(I,\left(0,0,...,0,\frac12\right)\right)\in E(n).$
From \cite[Lemma 3.2, 3.3]{KM} $s_{1}^{2},s_{2}^{2},...,s_{n}^{2}$ commute with each
other and generate a free abelian subgroup $\Z^n.$ In other words $M(A)$ is a flat manifold with holonomy group $Z_2^k$ of
diagonal type. Here $k$ is a number of non zero rows of a matrix $A$.

We have the following two lemmas.

\begin{lm}[\cite{KM}, Lemma 2.1]
	The cohomology ring $H^*(M(A),\mathbb Z_2)$ is generated  by degree one elements $x_1,\ldots,x_n$ as a graded ring with $n$ relations
	$$x_j^2=x_j\sum_{i=1}^na_{ij}x_i,$$
	for $j=1,\ldots,n$.
\end{lm}
\begin{lm}[\cite{KM}, Lemma 2.2]\label{lemma12}
	The real Bott manifold $M(A)$ is orientable if and only if the sum of entries is $0(\operatorname{mod}2)$ for each row of the matrix $A$.
\end{lm}

There are a few ways to decide whether there exists a Spin-structure
on an oriented flat manifold $M^n$. We start with
\begin{defi}[\cite{Fr}]
	An oriented flat manifold $M^n$ has a Spin-structure if and only if there exists a homomorphism $\epsilon\colon\Gamma\to\operatorname{Spin}(n)$ such that $\lambda_n\epsilon=p$, where $\lambda_n:\operatorname{Spin}(n)\to\operatorname{SO}(n)$ is the covering map.
\end{defi}
There is an equivalent condition for existence of Spin-structure. This is well known (\cite{Fr}) that the closed oriented
differential manifold $M$ has a Spin-structure if and only if the second Stiefel-Whitney class vanishes.

The $k$-th Stiefel-Whitney class \cite[ page 3, (2.1) ]{LS}  is given by the formula
\begin{equation}
w_k(M(A)) = (B(p))^{\ast}\sigma_{k}(y_1,y_2,...,y_{n})\in H^{k}(M(A);\Z_2) ,
\end{equation}
where $\sigma_k$ is the $k$-th elementary symmetric function, $B(p)$ is a map
induced by $p$ on the classification space and
\begin{equation}
y_i : = w_1(L_{i-1})\label{y}\end{equation}
for $i=2,3,\ldots,n$.
Hence,
\begin{equation}\label{sw1}
w_{2}(M(A)) = \sum_{1< i< j\leq n} y_{i}y_{j}\in H^{2}(M(A);\Z_2).
\end{equation}
\begin{defi} {(\cite{CMO}, page 4)}
	A binary square matrix $A$ is a Bott matrix if $A = PBP^{-1}$ for a permutation
	matrix $P$ and a strictly upper triangular binary matrix $B.$
\end{defi}

Our paper is a sequel of \cite{GS1}. There are given some conditions of the existence of Spin-structures.
\begin{theo}{\rm{(\cite{GS1}, page 1021)}}
	Let $A$ be a matrix of an orientable real Bott manifold $M(A)$.
	\begin{enumerate}
		\item Let $l\in\mathbb N$ be an odd number. If there exist $1\leq i<j\leq n$ and rows $A_{i,*}$, $A_{j,*}$ such that
		$$\sharp\{m:a_{i,m}=a_{j,m}=1\}=l$$
		and
		$$a_{ij}=0$$
		then $M(A)$ has no Spin-structure.
		\item If $a_{ij}=1$  and there exist $1\leq i<j\leq n$ and rows
		$$\begin{aligned}
		A_{i,*}&=(0,\ldots,0,a_{i,i_1},\ldots,a_{i,i_{2k}},0,\ldots,0),\\
		A_{j,*}&=(0,\ldots,0,a_{j,i_{2k+1}},\ldots,a_{j,i_{2k+2l}},0,\ldots,0)
		\end{aligned}$$
		such that $a_{i,i_1}=\ldots=a_{i,i_{2k}}=1$, $a_{i,m}=0$ for $m\not\in\{i_1,\ldots,i_{2k}\}$, $a_{j,i_{2k+1}}=\ldots=a_{j,i_{2k+2l}}=1$, $a_{j,r}=0$ for $r\not\in\{i_{2k+1},\ldots,i_{2k+2l}\}$ and $l$, $k$ are odd then $M(A)$ has no Spin-structure.
	\end{enumerate}
\end{theo}

In this paper we extend this theorem and we formulate necessary and sufficient conditions of the existence of a Spin-structure on real Bott manifolds. Here is our main result for Bott manifolds with holonomy group $Z_2^k$, $k$ even. Here is our main result
\begin{theo}\label{theolast}
	Let $A$ be a Bott matrix with $k$ non zero rows where $k$ is an even number. Then the real Bott manifold manifold $M(A)$ has a Spin-structure if and only if for all $1\leq i<j\leq n$  manifolds  $M(A_{ij})$ have a Spin-structure, where $A_{ij}$ is a matrix with $i-$ and $j-$ th nonzero rows.
\end{theo}

The structure of a paper is as follows. In Section 2 we give three lemmas. First of them gives a decomposition of the $n\times n-$integer matrix $A$ into $n\times n-$integer
matrices $A_{ij}$ with $i-$th and $j-$th nonzero rows. In Lemmas 2.2. and 2.3 we examine dependence of $y_i$ and $w_2$ of a real Bott manifold $M(A)$ on values $y_i^{jk}$ and $w_2(M(A_{jk}))$ of manifolds
$M(A_{jk})$. Then the proof of Theorem 1.2 will follow from Lemmas 2.2. and 2.3. Section 3 has a very technical character. In this section we shall give a complete characterization of the existence of the Spin-structure
on manifolds $M(A_{ij})$, $1\leq i<j\leq n$. Almost all statements in part 2 and 3 are illustrated by examples.

The author is grateful to Andrzej Szczepa\'{n}ski for his valuable, suggestions and help.

\section{Proof of the Main Theorem}

At the beginning we give formula for the decomposition of real Bott matrix $A$ into the sum of the real Bott matrices with two nonzero rows.

\begin{lm}
	Let $A$ be $n\times n-$Bott matrix and let $A_{ij}$, $1\leq i<j\leq n$, be $n\times n$-matrices with  $i-$th and $j-$th nonzero rows. Then, if $k$ is even, we have the following decomposition
	\begin{equation}A=\sum_{1\leq i<j\leq n}A_{ij}.\label{rozklad}\end{equation}
\end{lm}

\noindent
{\bf Proof.}
Let $A$ be $n\times n$-Bott matrix with $k$ nonzero rows, $k$ is an even number. Without loss of generality we can assume that nonzero rows have numbers from 1 to $k$. We shall consider the matrix $A$ as a sum of matrices $A_{ij}$, $1\leq i<j\leq n$. The number of matrices $A_{ij}$ is equal
$k\choose2$. For $1\leq i\leq k$ there are $(k-1)$-two elements subsets of $\{1,2,\ldots,k\}$ containing $i$.
Thus having summed matrices $A_{ij}$ we obtain
\begin{equation}
(k-1)\cdot A=\sum_{1\leq i<j\leq n}A_{ij}.\label{rozklad_1}
\end{equation}
Since $A$ is Bott matrix and $k$ is an even number we get
the formula (\ref{rozklad}).

\hskip 142mm $\Box$

\begin{ex}\label{ex1}
	Let
	$$A=\left[
	\begin{matrix}
	0&1&1&0&0&0\\
	0&0&1&1&0&0\\
	0&0&0&1&1&0\\
	0&0&0&0&1&1\\
	0&0&0&0&0&0\\
	0&0&0&0&0&0
	\end{matrix}\right]$$
	Thus $n=6$, $k=4$, so we have
	$$\begin{aligned}
	A=&
	\underbrace{\left[
		\begin{matrix}
		0&1&1&0&0&0\\
		0&0&1&1&0&0\\
		0&0&0&0&0&0\\
		0&0&0&0&0&0\\
		0&0&0&0&0&0\\
		0&0&0&0&0&0
		\end{matrix}\right]}_{A_{12}}
	+\underbrace{\left[
		\begin{matrix}
		0&1&1&0&0&0\\
		0&0&0&0&0&0\\
		0&0&0&1&1&0\\
		0&0&0&0&0&0\\
		0&0&0&0&0&0\\
		0&0&0&0&0&0
		\end{matrix}\right]}_{A_{13}}
	+\underbrace{\left[
		\begin{matrix}
		0&1&1&0&0&0\\
		0&0&0&0&0&0\\
		0&0&0&0&0&0\\
		0&0&0&0&1&1\\
		0&0&0&0&0&0\\
		0&0&0&0&0&0
		\end{matrix}\right]}_{A_{14}}\\
	&+\underbrace{\left[
		\begin{matrix}
		0&0&0&0&0&0\\
		0&0&1&1&0&0\\
		0&0&0&1&1&0\\
		0&0&0&0&0&0\\
		0&0&0&0&0&0\\
		0&0&0&0&0&0
		\end{matrix}\right]}_{A_{23}}
	+\underbrace{\left[
		\begin{matrix}
		0&0&0&0&0&0\\
		0&0&1&1&0&0\\
		0&0&0&0&0&0\\
		0&0&0&0&1&1\\
		0&0&0&0&0&0\\
		0&0&0&0&0&0
		\end{matrix}\right]}_{A_{24}}
	+\underbrace{\left[
		\begin{matrix}
		0&0&0&0&0&0\\
		0&0&0&0&0&0\\
		0&0&0&1&1&0\\
		0&0&0&0&1&1\\
		0&0&0&0&0&0\\
		0&0&0&0&0&0
		\end{matrix}\right]}_{A_{34}}
	\end{aligned}$$
\end{ex}
\vskip5mm
Before we start a proof of the main theorem we give an example.

\begin{ex}
	For the manifold $M(A)$ from Example \ref{ex1} we get
	$$
	y_2=x_1,\;
	y_3=x_1+x_2,\;
	y_4=x_2+x_3,\;
	y_5=x_3+x_4,\;
	y_6=x_4.$$
	Hence
	$$\omega_2(M(A))=x_1x_3+x_2x_4.$$
	
	We compute second Stiefel-Whitney classes for real Bott manifolds $M(A_{ij})$ from Example \ref{ex1}. For these purpose we put
	$y_l^{ij}=w_1(L_{l-1})$
	for manifolds $M(A_{ij})$ and we obtain
	$$\begin{array}{lllllll}
	y_2^{12}=x_1&y_2^{13}=x_1&y_2^{14}=x_1&y_2^{23}=0&y_2^{24}=0&y_2^{34}=0\\
	y_3^{12}=x_1+x_2&y_3^{13}=x_1&y_3^{14}=x_1&y_3^{23}=x_2&y_3^{24}=x_2&y_3^{34}=0\\
	y_4^{12}=x_2&y_4^{13}=x_3&y_4^{14}=0&y_4^{23}=x_2+x_3&y_4^{24}=x_2&y_4^{34}=x_3\\
	y_5^{12}=0&y_5^{13}=x_3&y_5^{14}=x_4&y_5^{23}=x_3&y_5^{24}=x_4&y_5^{34}=x_3+x_4\\
	y_6^{12}=0&y_6^{13}=0&y_6^{14}=x_4&y_6^{23}=0&y_6^{24}=x_4&y_6^{34}=x_4
	\end{array}$$
	With the above notation we get
	$$\begin{aligned}
	\sum_{1\leq i<j\leq k}y_2^{ij}&=3x_1=x_1\Rightarrow \sum_{1\leq i<j\leq k}y_2^{ij}=y_2,\\
	\sum_{1\leq i<j\leq k}y_3^{ij}&=3x_1+3x_2=x_1+x_2\Rightarrow \sum_{1\leq i<j\leq k}y_3^{ij}=y_3,\\
	\sum_{1\leq i<j\leq k}y_4^{ij}&=3x_2+3x_3=x_2+x_3\Rightarrow\sum_{1\leq i<j\leq k}y_4^{ij}=y_4,\\
	\sum_{1\leq i<j\leq k}y_5^{ij}&=3x_3+3x_4=x_3+x_4\Rightarrow\sum_{1\leq i<j\leq k}y_5^{ij}=y_5,\\
	\sum_{1\leq i<j\leq k}y_6^{ij}&=3x_4=x_4\Rightarrow\sum_{1\leq i<j\leq k}y_6^{ij}=y_6
	\end{aligned}$$
	and second Stiefel-Whitney classes for manifolds $M(A_{ij})$ are follows
	$$\begin{aligned}
	w_2(M(A_{12}))&=0,\\
	w_2(M(A_{13}))&=x_1x_3,\\
	w_2(M(A_{14}))&=0,\\
	w_2(M(A_{23}))&=0,\\
	w_2(M(A_{24}))&=0x_2x_4,\\
	w_2(M(A_{34}))&=0.
	\end{aligned}$$
	Hence
	$$
	\sum_{1=i<j=4}\omega_2(M(A_{ij})=x_1x_3+x_2x_4=\omega_2(M(A)).
	$$
	
\end{ex}

Following the method described in the above example we have lemmas.
\begin{lm}\label{lemat1}
	Let $A$ be a $n\times n$ Bott matrix with $k>3$ nonzero rows, $k$ is an even number.
	Then
	\begin{equation}y_l=\sum_{1\leq i<j\leq k}y_l^{ij},\label{wzor_y}\end{equation}
	where $y_l=\omega_1(L_{l-1}(M(A))$ and $y_l^{ij}=\omega_1(L_{l-1}(M(A_{ij}))$.
\end{lm}

\noindent
{\bf Proof.}
We have
$$y_l=w_1(L_{l-1})=\sum_{k=1}^{l-1}a_{kl}x_k=x\cdot A^l$$
where $x=[x_1,\ldots,x_n]$, $A=[a_{ij}]$, $A^l$ is the $l-$th column
of the matrix $A$ and $\cdot$ is multiplication of matrices. Let us multiply (\ref{rozklad})

on the left by $x$
$$\begin{aligned}
x\cdot A&=\sum_{1\leq i<j\leq k} x\cdot A_{ij}.\\
\end{aligned}$$
Since $yx\cdot A=[y_1,y_2,\ldots,y_n]$ and $x\cdot A^{ij}=[y^{ij}_1,y^{ij}_2,\ldots,y^{ij}_n]$, we get (\ref{wzor_y}).

\hskip 142mm $\Box$

\begin{lm}\label{lemat2}
	Let $A$ be $n\times n$ Bott matrix with $k-$nonzero rows, $k\geq4$, $k$ is an even number.
	Then
	$$w_2(M(A))=\sum_{1\leq i<j\leq k}w_2(M(A_{ij})).$$
\end{lm}

\noindent
{\bf Proof.}
From (\ref{sw1}) and (\ref{wzor_y})
$$\begin{aligned}
\omega_2(M(A))&=\sum_{l<r}y_ly_r\\
&=\sum_{l<r}\left[\left(\sum_{i<j}y_l^{ij}\right)\right]\left[\left(\sum_{i<j}y_r^{ij}\right)\right]
=\sum_{l<r}\left(\sum_{i<j}y_l^{ij}y_r^{ij}\right)\\
&=\sum_{i<j}\left(\sum_{l<r}y_l^{ij}y_r^{ij}\right)
=\sum_{i<j}\omega_2(M(A_{ij})).
\end{aligned}$$
\hskip 142mm $\Box$

From proofs of  Lemma \ref{lemat1} and Lemma \ref{lemat2} we obtain a proof of Main Theorem \ref{theolast}.

\noindent
{\bf Proof of Theorem \ref{theolast}}
Let us recall the manifold $M$ has a Spin-structure if and only if $w_2(M)=0$.
At the beginning let us assume, for each pair $1\leq i<j\leq n$, we have $w_2(M(A_{ij}))=0$. Then from Lemma \ref{lemat2} we have
$$w_2(M(A))=\sum_{1\leq i<j\leq k}w_2(M(A_{ij}))=0,$$
so the real Bott manifold $M(A)$ has a Spin-structure.

On the other hand, let the manifold $M(A)$ admits the Spin-structure, then
$$0=w_2(M(A))=\sum_{1\leq i<j\leq k}w_2(M(A_{ij})).$$
Second Stiefel-Whitney classes $M(A_{ij})$ are non negative so
$$\forall_{1\leq i<j\leq n}w_2(M(A_{ij}))=0.$$
\hskip 142mm $\Box$


\begin{rem}
	We do not know how to prove the main theorem for odd $k$.
	From the other side we are not sure if we can formulate it as a conjecture
	in this case.
\end{rem}

In the next section of our paper we concentrate on calculations of Spin-structure on manifolds $A_{ij}$.

\vskip 4mm
\section{Existence of Spin-structure on manifolds $M(A_{ij})$}

From now, let $A$ be a matrix of an orientable real Bott manifold $M(A)$ of dimension $n$ with two non-zero rows. From Lemma \ref{lemma12} we have that the number of entries 1, in each row, is an odd number and we have following three cases:
\newline {\bf CASE I.} There are no columns with double entries 1,
\newline{\bf CASE II.} The number of columns with double entries 1 is an odd number,
\newline {\bf CASE III.} The number of columns with double entries 1 is an even number,
\vskip5mm
\noindent
We give conditions for an existence of the Spin-structure on $M(A_{ij})$. In the further part of the paper we adopt the notation  $0_{p}=(\underbrace{0,\ldots,0}_{p\text{  - times}})$. From the definition, rows of number $i$ and $j$ correspond to generators $s_i,s_j$ which define a finite index abelian subgroup $H\subset\pi_1(M(A))$ (see \cite{GK}).

\begin{theo}\label{theo1}
	Let $A$ be a matrix of an orientable real Bott manifold $M(A)$ from the above case I. If there exist $1\leq i<j\leq n$ such that
	
	\noindent {\bf 1.}
	$$\begin{matrix}
	A_{i,\ast} &=(0_{i_1},a_{i,i_{1}+1},\ldots,a_{i,i_{1}+2k},0_{i_{2l}},0_{i_{p}})\\
	A_{j,\ast} &=(0_{i_1},0_{i_{2k}},a_{j,i_{1}+2k+1},\ldots,a_{j,i_{1}+2k+2l},0_{i_{p}}),\end{matrix}$$
	where $a_{i,i_{1}+1} = \ldots = a_{i,i_{1}+2k} = 1, a_{i,m} = 0$ for $m\notin\{i_1,\ldots,i_{1}+2k\}$,
	$a_{j,i_{1}+2k+1}= \ldots = a_{j,i_{1}+2k+2l} = 1, a_{j,r} = 0$ for
	$r\notin\{i_{1}+2k+1,\ldots,i_{1}+2k+2l\}$.
	\newline Then $M(A)$ admits the Spin-structure if and only if either $l$ is an even number or $l$ is an odd number and
	$j\notin\{i_1+1,\ldots,i_{1}+2k\}$.
	
	\noindent
	{\bf 2.}
	$$\begin{matrix}
	A_{i,\ast} &=(0_{i_{1}},0_{i_{2k}},a_{i,i_{2k}+1},\ldots,a_{i,i_{2k}+2l},0_{i_{p}})\\
	A_{j,\ast} &=(0_{i_{1}},a_{j,i_{1}+1},\ldots,a_{j,i_{1}+2k},0_{i_{2l}},0_{i_{p}}),\end{matrix}$$
	where
	$a_{j,i_{1}+1} = \ldots = a_{j,i_{1}+2k} = 1, a_{j,m} = 0$ for $m\notin\{i_1,\ldots,i_{1}+2k\}$,
	$a_{i,i_{2k}+1} = \ldots = a_{i,i_{2k}+2l} = 1, a_{i,r} = 0$ for
	$r\notin\{i_{2k}+1,\ldots,i_{2k}+2l\}$,
	then $M(A)$ has the Spin-structure.
\end{theo}

\noindent
{\bf Proof.}
{\bf 1.}
From (\ref{y}) we have
$$\begin{aligned}
y_{i_1+1}&=\ldots=y_{i_{1}+2k}=x_i,\\
y_{i_{1}+2k+1}&=\ldots=y_{i_{1}+2k+2l}=x_j.\end{aligned}$$
Using (\ref{sw1}) and
$x_i^2=x_i\sum_{j=1}^na_{ji}x_j$
we get
$$\begin{aligned}w_2(M(A))&=k(2k-1)x_i^2+4klx_ix_j+l(2l-1)x_j^2\\
&=k(2k-1)x_i^2+l(2l-1)x_j^2=l(2l-1)x_j^2=lx_j^2.
\end{aligned}$$
Summing up, we have to consider the following cases
\begin{enumerate}
	\item if $l=2b$, then $w_2(M(A))=2bx_j^2=0$. Hence $M(A)$ has a Spin-structure,
	\item if $l=2b+1$, then
	$$\begin{aligned}
	w_2(M(A))&=(2b+1)x_j^2=x_j^2\\
	&=\begin{cases}0,&\text{if }j\notin\{i_1+1,\ldots,i_{1}+2k\},M(A)\text{ has a Spin-structure,}\\x_ix_j,&\text{if }j\in\{i_1+1,\ldots,i_{1}+2k\},M(A)\text{ has no Spin-structure}.\end{cases}
	\end{aligned}$$
\end{enumerate}
{\bf 2.}
From {\rm(}\ref{y}{\rm)}
$$\begin{aligned}
y_{i_1}+1&=\ldots=y_{i_{1}+2k}=x_j\\
y_{i_{1}+2k+1}&=\ldots=y_{i_{1}+2k+2l}=x_i.
\end{aligned}$$
Moreover, from (\ref{sw1}) and since  $i_1>j>i$
$$\begin{aligned}
w_2(M(A))&=k(2k-1)x_j^2+4klx_ix_j+l(2l-1)x_i^2\\
&=k(2k-1)\underbrace{x_j^2}_{=0}+l(2l-1)\underbrace{x_i^2}_{=0}=0.\end{aligned}$$
Hence $M(A)$ has the Spin-structure.

\hskip 142mm $\Box$

\begin{theo}\label{theo2} Let $A$ be a matrix of an orientable real Bott manifold $M(A)$ from the above case II. If there exist $1\leq i<j\leq n$ such that
	\newline{\bf 1.}
	$$\begin{footnotesize}\begin{matrix}
	A_{i,\ast} = (0_{i_1},a_{i,i_{1}+1},\ldots,a_{i,i_{1}+2k},a_{i,i_{1}+2k+1},\ldots,a_{i,i_{1}+2k+2l},0_{i_{2m}},0_{i_p})\\
	A_{j,\ast} = (0_{i_1},0_{i_{2k}},a_{j,i_{1}+2k+1},\ldots,a_{j,i_{1}+2k+2l},a_{j,i_{1}+2k+2l+1},\ldots,a_{j,i_{1}+2k+2l+2m},0_{i_p}),\end{matrix}\end{footnotesize}$$
	where
	$a_{i,i_{1}+1} = \ldots = a_{i,i_{1}+2k} =a_{i,i_{1}+2k+1}=\ldots=a_{i,i_{1}+2k+2l}=1, a_{i,r} = 0$ for $r\notin\{i_1+1,\ldots,i_{1}+2k+2l\}$,
	$a_{j,i_{1}+2k+1} = \ldots= a_{j,i_{1}+2k+2l+2m} = 1, a_{j,s} = 0$ for
	$s\notin\{i_{1}+2k+1,\ldots,i_{1}+2k+2l+2m\}$.
	
	Then $M(A)$ has the Spin-structure if and only if either $l$ and $m$ are number of the same parity or $l$ and $m$ are number of different parity and
	$j\notin\{i_1+1,\ldots,i_{1}+2k\}$.
	
	\noindent {\bf 2.}
	$$\begin{footnotesize}\begin{matrix} A_{i,\ast} &=(0_{i_1},0_{i_{1}+2k},a_{i,i_{1}+2k+1},\ldots,a_{i,i_{1}+2k+2l},a_{i,i_{1}+2k+2l+1},\ldots,a_{i,i_{1}+2k+2l+2m},0_{i_p}) ,\\
	A_{j,\ast} &=(0_{i_1},a_{j,i_{1}+1},\ldots,a_{j,i_{1}+2k},a_{j,i_{1}+2k+1},\ldots,a_{j,i_{1}+2k+2l},0_{i_{2m}},0_{i_p}) \end{matrix}\end{footnotesize}$$
	where
	$a_{j,i_{1}+1}= \ldots = a_{j,i_{1}+2k}=a_{j,i_{1}+2k+1}=\ldots=a_{j,i_{1}+2k+2l} = 1, a_{j,m} = 0$ for $m\notin\{i_1+1,\ldots,i_{1}+2k+2l\}$,
	$a_{i,i_{1}+2k+1}= \ldots = a_{i,i_{1}+2k+2l}=a_{i,i_{1}+2k+2l+1}=\ldots=a_{i,i_{1}+2k+2l+2m} = 1, a_{i,r} = 0$ for
	$r\notin\{i_{1}+2k+1,\ldots,i_{1}+2k+2l+2m\}$, then $M(A)$ has the Spin-structure
\end{theo}

\noindent{\bf Proof.}
{\bf 1.} From (\ref{y}) we have
$$\begin{aligned}
y_{i_1+1}&=\ldots=y_{i_{1}+2k}=x_i,\\
y_{i_{1}+2k+1}&=\ldots=y_{i_{1}+2k+2l}=x_i+x_j\\
y_{i_{1}+2k+2l+1}&=\ldots=y_{i_{1}+2k+2l+2m}=x_j.\end{aligned}$$
From (\ref{sw1}) and $x_i^2=x_i\sum_{j=1}^na_{ji}x_j$
we get
$$\begin{aligned}w_2(M(A))&=k(2k-1)x_i^2+4klx_i(x_i+x_j)+l(2l-1)(x_i+x_j)^2+m(2m-1)x_j^2\\
&=l(2l-1)x_j^2+m(2m-1)x_j^2=(l+m)x_j^2.
\end{aligned}$$
We have to consider the following cases:
\begin{enumerate}
	\item If $l+m$ is an even number then $w_2(M(A))=0.$ Hence $M(A)$ has a Spin-structure.
	\item If $l+m$ is an odd number then
	$$\begin{aligned}
	w_2(M(A))&=x_j^2\\
	&=\begin{cases}0,&\text{if }j\notin\{i_1+1,\ldots,i_{1}+2k\},M(A)\text{ has a Spin-structure}\\x_ix_j,&\text{if }j\in\{i_1+1,\ldots,i_{1}+2k\},M(A)\text{ has no Spin-structure}.\end{cases}
	\end{aligned}$$
\end{enumerate}
{\bf 2.} Using (\ref{y}) we get
$$\begin{aligned}
y_{i_1+1}&=\ldots=y_{i_{1}+1}=x_j\\
y_{i_{1}+2k+1}&=\ldots=y_{i_{1}+2k+2l}=x_i+x_j\\
y_{i_{1}+2k+2l+1}&=\ldots=y_{i_{1}+2k+2l+2m}=x_i.
\end{aligned}$$
Moreover, from {\rm(}\ref{sw1}{\rm)} and  since $i_1>j>i$
$$\begin{aligned}
&w_2(M(A))=k(2k-1)x_j^2+l(2l-1)x_i^2+4klx_j(x_i+x_j)+4kmx_ix_j\\
&+4lmx_i(x_i+x_j)+l(2l-1)(x_i+x_j)^2+m(2m-1)x_i^2\\
&=k(2k-1)\underbrace{x_j^2}_{=0}+l(2l-1)\underbrace{x_i^2}_{=0}+l(2l-1)\underbrace{x_j^2}_{=0}+m(2m-1)\underbrace{x_i^2}_{=0}=0.
\end{aligned}$$
Hence $M(A)$ has a Spin-structure.

\hskip 142mm $\Box$

\begin{theo}\label{theo3} Let $A$ be a matrix of an orientable real Bott manifold $M(A)$ from the above case III. If there exist $1\leq i<j\leq n$ such that
	\newline {\bf 1.}
	$$\begin{footnotesize}\begin{matrix} A_{i,\ast} &= (0_{i_1},a_{i,i_{1}+1},\ldots,a_{i,i_{1}+2k+1},a_{i,i_{1}+2k+2},\ldots,a_{i,i_{1}+2k+2l+2},0_{i_{2m+1}},0_{i_p})\\
	A_{j,\ast} &= (0_{i_1},0_{i_{2k+1}},a_{j,i_{2k+2}},\ldots,a_{j,i_{1}+2k+2l+2},a_{j,i_{1}+2k+2l+3},\ldots,a_{j,i_{1}+2k+2l+2m+3},0_{i_p}),\end{matrix}\end{footnotesize}$$
	where
	$a_{i,i_{1}+1} = \ldots = a_{i,i_{1}+2k} =\ldots=a_{i,i_{1}+2k+2l+2}=1, a_{i,r} = 0$ for $r\notin\{i_1+1,\ldots,i_{1}+2k+2l+2\}$,
	$a_{j,i_{1}+2k+2} = \ldots= a_{j,i_{1}+2k+2l+2m+3} = 1, a_{j,s} = 0$ for
	$s\notin\{i_{1}+2k+2,\ldots,i_{1}+2k+2l+2m+3\}$.
	Then $M(A)$ admits the Spin-structure if and only $l$ and $m$ are number of the same parity and
	$j\in\{i_1+1,\ldots,i_{1}+2k+2\}$ .
	
	{\bf 2.}
	$$\begin{footnotesize}\begin{matrix} A_{i,\ast} &=(0,_{i_1},0_{i_{2l+1}},a_{i,i_{1}+2k+2},\ldots,a_{i,i_{1}+2k+2l+2},a_{i,i_{1}+2k+2l+3},\ldots,a_{i,i_{1}+2k+2l+2m+3},0_{i_p})\\
	A_{j,\ast} &=(0_{i_1},a_{j,i_{1}+1},\ldots,a_{j,i_{1}+2k+1},a_{j,i_{1}+2k+2},\ldots,a_{j,i_{1}+2k+2l+2},0_{i_{2m}},0_{i_p}) \end{matrix}\end{footnotesize}$$
	where
	$a_{j,i_{1}+1}= \ldots = a_{j,i_{1}+2k}=a_{j,i_{1}+2k+1}=\ldots=a_{j,i_{1}+2k+2l+2} = 1, a_{j,m} = 0$ for $m\notin\{i_1+1,\ldots,i_{1}+2k+2l+2\}$,
	$a_{i,i_{1}+2k+2}= \ldots = a_{i,i_{1}+2k+2l+2}=a_{i,i_{1}+2k+2l+3}=\ldots=a_{i,i_{1}+2k+2l+2m+3} = 1, a_{i,r} = 0$ for
	$r\notin\{i_{1}+2k+2,\ldots,i_{1}+2k+2l+2m+3\}$.
	Then $M(A)$ has no Spin-structure.
\end{theo}

\noindent
{\bf Proof.} {\bf 1.}
From (\ref{y})
$$\begin{aligned}
y_{i_1+1}&=\ldots=y_{i_{1}+2k+1}=x_i,\\
y_{i_{1}+2k+2}&=\ldots=y_{i_{1}+2k+2l+2}=x_i+x_j\\
y_{i_{1}+2k+2l+3}&=\ldots=y_{i_{1}+2k+2l+2m+3}=x_j.\end{aligned}$$
From (\ref{sw1}) and $x_i^2=x_i\sum_{j=1}^na_{ji}x_j$
we obtain
$$\begin{aligned}w_2(M(A))&=k(2k+1)x_i^2+(2k+1)(2l+1)x_i(x_i+x_j)+(2k+1)(2m+1)x_ix_j\\
&+l(2l+1)(x_i+x_j)^2+(2l+1)(2m+1)x_j(x_i+x_j)+m(2m+1)x_j^2\\
&=(l+m+1)x_j^2+(2l+1)(2m+1)x_ix_j=(l+m+1)x_j^2+x_ix_j.
\end{aligned}$$
Now, if $l$ and $m$ are number of the same parity we have
$$\begin{aligned}&w_2(M(A))=x_ix_j+x_j^2\\&=
\begin{cases}x_ix_j,&\text{ if } j\notin\{i_1+1,\ldots,i_{1}+2k+2\}, M(A)\text { has no Spin-structure},\\
0,&\text{ if } j\in\{i_1+1,\ldots,i_{1}+2k+2\}, M(A)\text  { has a Spin-structure}.\end{cases}
\end{aligned}$$

{\bf 2.} From {\rm(}\ref{y}{\rm)}
$$\begin{aligned}
y_{i_1+1}&=\ldots=y_{i_{1}+2k+1}=x_j\\
y_{i_{1}+2k+2}&=\ldots=y_{i_{1}+2k+2l+2}=x_i+x_j\\
y_{i_{1}+2k+2l+3}&=\ldots=y_{i_{1}+2k+2l+2m+3}=x_i.
\end{aligned}$$
From (\ref{sw1}) and since $i_1>j>i$ we get
$$\begin{aligned}
w_2(M(A))&=k(2k+1)x_j^2+m(2m+1)x_i^2+(2k+1)(2l+1)x_j(x_i+x_j)\\
&+(2k+1)(2m+1)x_ix_j+l(2l+1)(x_i+x_j)^2\\
&+(2l+1)(2m+1)x_i(x_i+x_j)+m(2m-1)x_i^2\\
&=k(2k+1)\underbrace{x_j^2}_{=0}+l(2l+1)\underbrace{(x_i+x_j)^2}_{=0}+m(2m+1)\underbrace{x_i^2}_{=0}\\
&+x_j(x_i+x_j)+x_ix_j+x_i(x_i+x_j)=x_ix_j\ne0,
\end{aligned}$$
so $M(A)$ has no Spin-structure.

\hskip 142mm $\Box$

Now, we give examples which illustrate Theorems \ref{theo1} - \ref{theo3}.
\begin{ex}
	
	{\bf 1.}	Let
	$$A=\left[\begin{matrix}
	0&0&0&0&0&0&0&0\\
	0&0&1&1&0&0&0&0\\
	0&0&0&0&1&1&1&1\\
	0&0&0&0&0&0&0&0\\
	0&0&0&0&0&0&0&0\\
	0&0&0&0&0&0&0&0\\
	0&0&0&0&0&0&0&0\\
	0&0&0&0&0&0&0&0
	\end{matrix}\right].$$
	Here $2l=4\Rightarrow l=2$. Hence from Theorem \ref{theo1}, part 1.1, manifold $M(A)$ has the Spin-structure.
	\newline{\bf 2.}
	$$A=\left[\begin{matrix}
	0&0&0&0&0&0\\
	0&0&1&1&0&0\\
	0&0&0&0&1&1\\
	0&0&0&0&0&0\\
	0&0&0&0&0&0\\
	0&0&0&0&0&0
	\end{matrix}\right].$$
	Here
	$
	l=1,\{i_1,i_2,\ldots,i_n\}=\{3,4\},j=3\in\{3,4\}.$
	Hence, from Theorem \ref{theo1}, part 1.2, the real Bott manifold $M(A)$ has no Spin-structure.
	\newline{\bf 3.}
	$$A=\left[\begin{matrix}
	0&1&1&1&1&0&0&0&0\\
	0&0&0&1&1&1&1&1&1\\
	0&0&0&0&0&0&0&0&0\\
	0&0&0&0&0&0&0&0&0\\
	0&0&0&0&0&0&0&0&0\\
	0&0&0&0&0&0&0&0&0\\
	0&0&0&0&0&0&0&0&0\\
	0&0&0&0&0&0&0&0&0\\
	0&0&0&0&0&0&0&0&0
	\end{matrix}\right].$$
	From Theorem \ref{theo2}, part 1.4 and since $l=1, m=2, \{i_1,\ldots,i_{2k}\}=\{2,3\}, j=2\in\{2,3\}$
	the real Bott manifold has no Spin-structure.
	\newline{\bf 4.}
	$$A=\left[\begin{array}{ccccccccccccc}
	0&0&0&0&0&0&0&0&0&0&0&0&0\\
	0&0&1&1&1&1&1&1&0&0&0&0&0\\
	0&0&0&0&0&1&1&1&1&1&1&1&1\\
	0&0&0&0&0&0&0&0&0&0&0&0&0\\0&0&0&0&0&0&0&0&0&0&0&0&0\\
	0&0&0&0&0&0&0&0&0&0&0&0&0\\0&0&0&0&0&0&0&0&0&0&0&0&0\\
	0&0&0&0&0&0&0&0&0&0&0&0&0\\0&0&0&0&0&0&0&0&0&0&0&0&0\\
	0&0&0&0&0&0&0&0&0&0&0&0&0\\0&0&0&0&0&0&0&0&0&0&0&0&0\\
	0&0&0&0&0&0&0&0&0&0&0&0&0\\0&0&0&0&0&0&0&0&0&0&0&0&0
	\end{array}\right].$$
	In this case $l=1, m=2, $
	and from Theorem \ref{theo3} we have that $M(A)$ has no Spin-structure.
\end{ex}

\vskip 2mm
\noindent
Maria Curie-Sk{\l}odowska University,\\
Institute of Mathematics\\
pl. Marii Curie-Sk{\l}odowskiej 1\\
20-031 Lublin, Poland\\
E-mail: anna.gasior@poczta.umcs.lublin.pl

\end{document}